\documentclass[12pt]{article}

\title{Positive speed for high-degree automaton groups}
\date{February 2011}

\author{Gideon Amir  \and B\'alint Vir\'ag}

\usepackage{amsmath,amsthm,amssymb,amsfonts,enumerate,graphicx}

\usepackage[longnamesfirst]{natbib}
\bibpunct[ ]{(}{)}{,}{a}{}{,}

\theoremstyle{plain}
    \newtheorem{theorem}{Theorem}
    \newtheorem{fact}[theorem]{Fact}
    \newtheorem{lemma}[theorem]{Lemma}
    \newtheorem{proposition}[theorem]{Proposition}

    \newtheorem{conjecture}[theorem]{Conjecture}
    \newtheorem{question}[theorem]{Question}

\theoremstyle{definition} 
    \newtheorem{definition}[theorem]{Definition}

\newcommand{\one}{{\mathbf 1}}
\newcommand{\M}{\mathfrak M}
\renewcommand{\L}{\mathfrak L}
\newcommand{\tree }{\mathbb T}

\newcommand{\energy}{\mathcal E}

\newcommand{\G}{\mathcal G}

\newcommand{\Aut}{\operatorname{Aut}}
\newcommand{\Sym}{\operatorname{Sym}}
\newcommand{\stab}{\operatorname{stab}}
\newcommand{\llangle}{\langle\hspace{-0.2em}\langle}
\newcommand{\rrangle}{\rangle\hspace{-0.2em}\rangle}
\newcommand{\ang}[1]{{\llangle #1 \rrangle}}

\renewcommand{\bar}{\overline}
\newcommand{\lamp}{\operatorname{lamp}}
\newcommand{\germ}{\operatorname{germ}}

\newcommand{\res}{\mathfrak R}
\newcommand{\rr}{\overline r}
\newcommand{\oo}{\overline o}

\begin{document}
\maketitle
\begin{abstract}
Mother groups are the basic building blocks for polynomial
automaton groups.

We show that, in contrast with mother groups of degree 0 or
1, any bounded, symmetric, generating random walk on the
mother groups of degree at least 3 has positive speed.

The proof is based on an analysis of resistance in fractal
mother graphs. We give upper bounds on resistances in these
graphs, and show that infinite versions are transient.
\end{abstract}
\begin{center}
\includegraphics*[width=400pt]{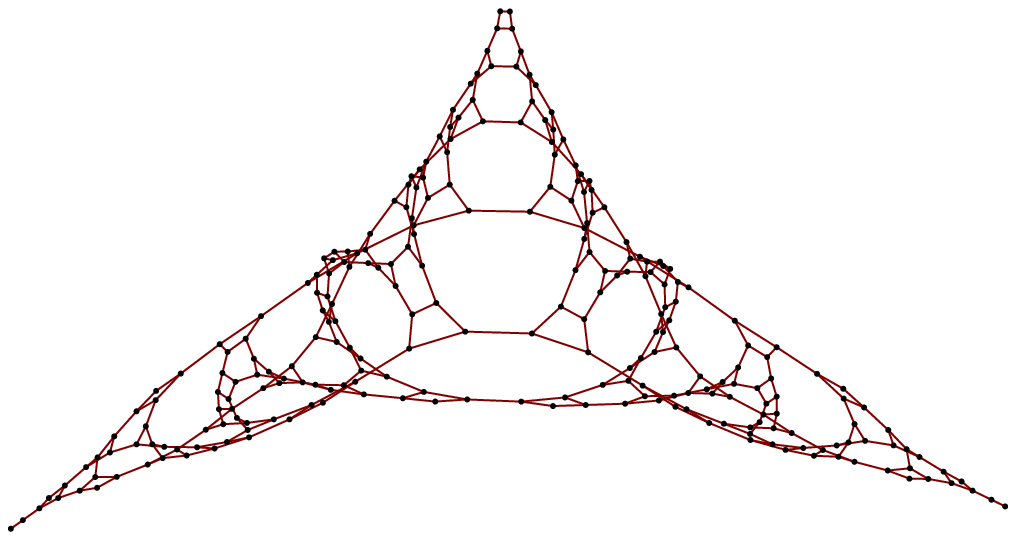}\\
The fractal mother graph $\G(1,2,8)$
\end{center}

\section{Introduction}

The most interesting automaton groups, called polynomial
automata, are classified according to their degree.
\cite{Sidki00} showed all these groups are small in the
sense that they do not contain nonabelian free subgroups.
He asked whether all these groups were amenable, a question
that is still wide open.

In \cite{BKN} (for degree 0) and \cite{AAV} (for all
degrees) it was shown that every polynomial automaton group
is contained in a mother group of the same degree. The
mother groups $\M_{d,n}$, defined precisely in Section
\ref{s:mother}, are a simple set of groups indexed by their
degree $d$ and the size of the alphabet $m$. Our main
theorem shows that for $d\ge 4$ or $d\ge 3, m\ge 3$ the
mother groups are large in the following sense.

\begin{theorem}\label{t:main} Every bounded, generating random walk on the mother groups with
$d\ge 4$ or $d\ge 3, m\ge 3$ has positive speed.
\end{theorem}

This is proved through Proposition \ref{p:mainresbound} and
Theorem \ref{t:liouville-tran}. Theorem \ref{t:main} is in
contrast with $d=0,1$ where the mother groups support
bounded, generating random walks with zero speed.
(\cite{BKN}, \cite{AAV}.)

The proof of this theorem is based on the analysis of
Schreier graphs of the natural action of automata groups.
These graphs are discrete versions of fractals, and indeed,
many classical fractals can be represented this way.
Examples include the Sierpinski gasket, the long-range
graph of \cite{BH05}, the Basilica (see \cite{gz02b} and
\cite{BV05}) and other Julia sets. The Schreier graphs
$\G(d,m,n)$ for the mother groups, called {\bf fractal
mother graphs} can be described up to uniform
quasi-isometry, as follows (the precise definition is given
in Section \ref{s:mother}). The vertex set is
$\{0,\ldots,m-1\}^n$, and two vertices are connected by an
edge if they differ in a single digit, and this digit is
followed by at most $d+1$ nonzero digits. Figures of
various fractal mother groups can be found throughout the
paper.
\begin{theorem}
The resistance between any two vertices in $\G(d,m,n)$ is
bounded above by
$$\begin{array}{ll}
c\, \left(\frac{m}{m-1}\right)^s &\text{for }  d=0,\\
c\,s^{3-d-\log_m (m-1)} &\text{for }  d=1 \text{ or }  d=2,\\
\,c\log^2 n &\text{for } d=3,m=2,\\
 c &\text{for } d\ge 4  \text{ or }d\ge 3, m \ge 3.
     \end{array}
$$
where  $c$ is a constant depending on $d,m$ only.
\end{theorem}

This will be shown in Propositions \ref{p:degree0bound} and
\ref{p:mainresbound}. It is interesting to note that the
graphs $\G(0,2,n)$ are paths, and $\G(0,3,n)$ are close
relatives of the Sierpinski gasket, although different
enough that the resistances are off by a power. See
\cite{NT08} for an overview of analysis on such fractals.

In an upcoming paper we plan to show a logarithmic lower
bound for $d=2$ and a polynomial lower bound for $d=1$.

\begin{figure}[t!]
\begin{center}
\includegraphics*[width=400pt]{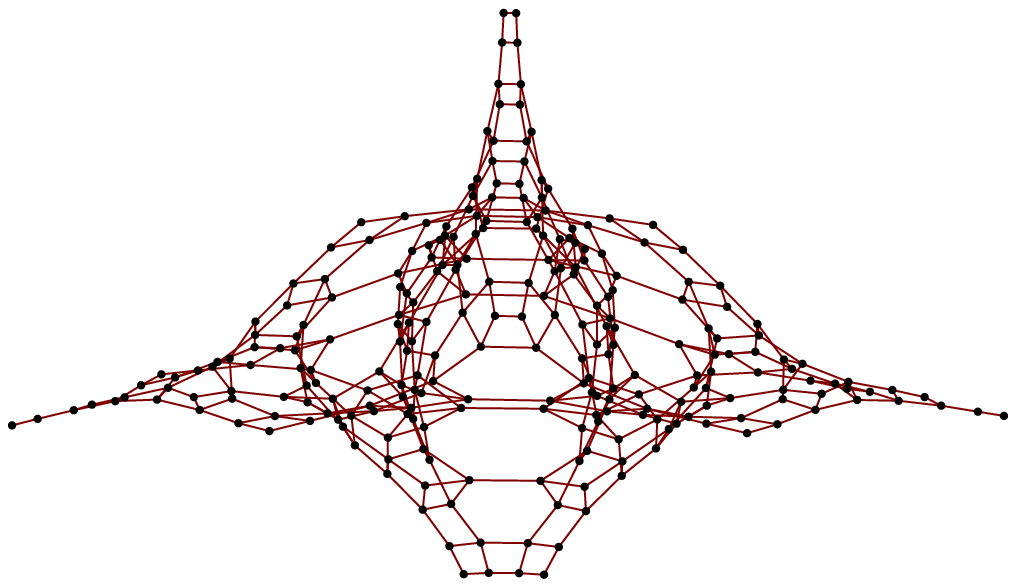}\\
The fractal mother graph $\G(2,2,8)$
\end{center}
\end{figure}
The graphs $\G(d,m,n)$ have infinite versions, defined as
above with $n=\infty$. These graphs have uncountably many
connected components; the notation $\G(d,m,\infty)$ refers
to the component of the vertex $\cdots 000$. In this case,
we show that for $d\ge 4$ or $d\ge 3, m\ge 3$ every
connected component is transient.

The second step of the proof is based on the ideas of
\cite{Erschler04}. In our setting we have the following.
\begin{theorem}\label{t:liouville-tran} If $\G(d,m,\infty)$ is transient, then
every bounded generating random walk on the group
$\M_{d,m}$ has positive speed.
\end{theorem}
In fact, following the proof of \cite{Erschler04}, we show
that the group has nonconstant bounded harmonic functions,
which is equivalent to positive speed, see \cite{KV},
\cite{varoplongrange}, \cite{vershik00} and also the
lecture notes \cite{pete}.

Finally, we would like to state some open questions and
conjectures. The one we wish could be a theorem in this
paper is
\begin{figure}[b!]
\begin{center}
\includegraphics*[width=350pt]{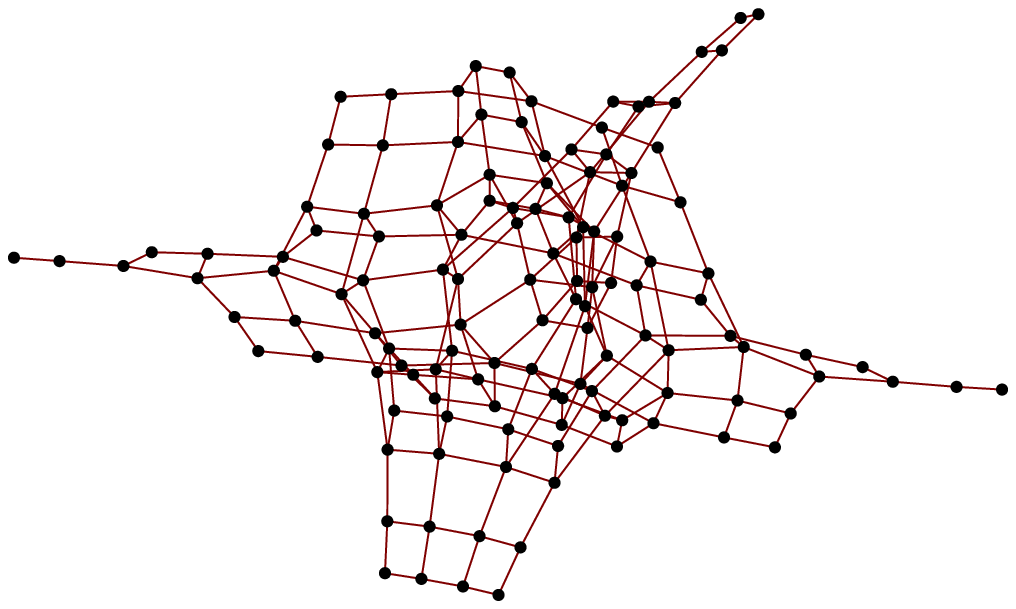}\\
The fractal mother graph $\G(3,2,7)$
\end{center}
\end{figure}
\begin{conjecture}
The graph $\G(3,2,\infty)$ is transient.
\end{conjecture}

This would take care of the missing case $d=3$, $m=2$.

Our bounds for $d=0$, $d\ge 4$ and $d\ge 3, m\ge 3$ are
sharp. We believe

\begin{conjecture}
The bound for $d=1$ is sharp.
\end{conjecture}

In contrast, it seems that we have the following.

\begin{conjecture}
The maximal resistance in $\G(2,m,n)$ grows slower than any
power of $n$.
\end{conjecture}

In a future paper, we will show that $\G(2,m,\infty)$ is
recurrent. This, however, does not decide

\begin{question}
Do the degree 2 mother groups support generating random
walks with positive speed? Do they support generating walks
with zero speed?
\end{question}

We have also observed, but have not been able to prove, the
following. In the graph $\G(d,m,n)$ call the vertex
$00\ldots 00$ the {\bf root} and the vertex $10\ldots 00$
the {\bf antiroot}.
\begin{conjecture}\label{maxres}
The root and the antiroot maximize the resistance in
$\G(d,m,n)$.
\end{conjecture}
Showing this would simplify our present proof.

\section{Automata and their groups}\label{s:automata}

\paragraph{Basic definitions.} Finite automata are the simplest interesting
model of computing. The space of words in alphabet
$\{0,\dots,m-1\}$ has a natural tree structure, with
$\{wx\}_{x<m}$ being the children of the finite word $w$,
and the empty word $\emptyset$ being the root. Let $\tree
_m$ denote this tree.

A {\bf finite automaton} on $m$ symbols is a finite set of states
$A$ with a map $A \to A^{m} \times \Sym(m)$ sending $a \mapsto
(a_0,\ldots, a_{m-1}, \sigma_a)$. We will use the notation
\[
a = \llangle  a_0,\ldots, a_{m-1} \rrangle  \sigma_a.
\]

An automaton acts on words in alphabet $\{0,\ldots,{m-1}\}$
sequentially. When the automaton is in a state $a$ and reads a
letter $x$, it outputs $x.\sigma_a$ and moves to state $a_x$. From
this state the automaton acts on the rest of the word.
Formally, for a word $xw$ (starting with a letter $x$) we have
the recursive definition
\begin{equation}\label{eq:group_action}
  (xw).a = (x.\sigma_a)(w.a_x).
\end{equation}
The first $k$ symbols of the output are determined by the first
$k$ symbols read. Note that the action is defined for for both
finite and infinite words, and that the action on infinite words
determines the action on finite words and vice versa. It follows
that each element $a\in A$ is an automorphism of $\tree _m$. The
{\bf automaton group} corresponding to an automaton $A$ is the
subgroup of $\Aut(\tree_m)$ generated by $A$.

The action \eqref{eq:group_action} corresponds to the following
{\bf multiplication rule}:
\[
\ang{a_0,\dots,a_d} \sigma \ang{b_0,\dots,b_d} \tau =
\ang{a_0b_{0.\sigma},\ldots,a_{m-1}b_{(m-1).\sigma}} \sigma
\tau.
\]
This multiplication rule can be used to define automaton
groups without any reference to automorphisms of the tree.

We use the conjugation notation $a^b = b^{-1}ab.$

The notion of first-level sections can be generalized to any
level. If $v\in \tree _m$ is a finite word and $g\in Aut(\tree
_m)$, then there is a word $v'$ of equal length to $v$ and an
automorphism $g'\in\Aut(\tree _m)$ such that $vw.g=v'(w.g')$, for
every word $w$. This $g'$ is called the {\bf section} of $g$ at
$v$. Informally, $g'$ is the action of $g$ on the subtree above
the vertex $v$. The section of $g$ at $v$ is denoted $g(v)$.

The group $\Aut(\tree_m)$ also acts on infinite geodesics
starting from the root, called {\bf rays}; these correspond
to one-sided infinite words in the alphabet. The set of
rays is called the boundary $\partial \tree_m$ of
$\tree_m$. The {\bf zero ray} is the geodesic $\ldots
0000$.

\paragraph{Activity growth of automaton groups.}

For any state $a\in A$, the number of length-$n$ words $v$ such
that the section $a(v)$ is not the identity satisfies a certain
linear recursion. Thus this number grows either polynomially with
some degree $d$ or exponentially. We define the {\bf degree of
activity growth} (in short, degree) of $a$ to be $d$ or $\infty$,
respectively. The {\bf degree} of an automaton group is the
maximal degree of any of its generators. Automaton groups are said
to have {\bf bounded, linear, polynomial} or {\bf exponential}
activity growth when their degree is $0$, $1$, finite or infinite,
respectively.

\section{Mother groups}\label{s:mother}

The {\bf mother group}, denoted $\M_{d,m}$, is defined as
the automaton group generated by the following states
\begin{align*}
  \alpha_{k,\sigma} &= \llangle  \alpha_{k,\sigma},
  \alpha_{k-1,\sigma},1,\ldots,1\rrangle , \qquad 0\le k \le d,  \\
  \alpha_{-1,\sigma} &= \sigma  \\
  \beta_{k,\rho} &= \llangle  \beta_{k,\rho},
  \beta_{k-1,\rho},1,\ldots,1\rrangle ,\qquad 1\le k \le d,\\
  \beta_{0,\rho} &= \llangle
  \beta_{0,\rho},1,\ldots,1\rrangle \rho.
\end{align*}
where $\sigma,\rho \in \Sym(m)$ are arbitrary, subject to
$0.\rho = 0$. The number of states in $\M_{d,m}$ as defined
here is $m!(d+1)+(m-1)!d$. The same group can be generated
with fewer states by taking $\sigma,\rho$ only in a
minimal, ($2$-element) set of generators of $\Sym(m)$ and
$\stab(0)\subset\Sym(m)$ respectively. This would give a
generating set of size $4d+2$.

The actions of $\alpha_{k,\sigma}$ and $\beta_{k,\rho}$ on a word
have simple descriptions. Both read the word and make no changes
up to the $k+1$'th non-zero letter.
\begin{itemize}
\item If the first $k+1$ nonzero letters in a word are all 1,
    then $\alpha_{k,\sigma}$ permutes the next letter by
    $\sigma$. Otherwise it does nothing.
\item If the first $k$ nonzero letters in a word are all 1,
    then $\beta_{k,\rho}$ permutes the next nonzero letter by
    $\rho$. Otherwise it does nothing.
\end{itemize}
Thus both affect only the $k+1$'st non-zero letter and the letter
immediately following it.

\begin{theorem}[Mother groups contain all, \cite{AAV}]\label{T:universal_mother}
  Every degree-$d$ automaton group is isomorphic to a subgroup of
  $\M_{d,m}$ for some $m$.
\end{theorem}

Note that $m$ is generally not the same as the degree of the tree
on which the automaton acts.

\subsection*{Level subgroups}

Observe that the group of automorphisms of the first two
levels of $\tree _m$ fixing $0$ and its children is
isomorphic to $\Sym(m) \wr \Sym(m-1)$. (We will interpret elements in $\Sym(m-1)$ as acting on $\{1,\ldots,m-1\}$.)

For each $\sigma\in \Sym(m) \wr \Sym(m-1)$ and each word $w$ in the symbols
$\{1,\dots,m-1\}$ let $\lambda_{w,\sigma}$ denote the element of $\Aut(\tree _m)$
acting as follows: If the first $|w|$ nonzero letters agree with $w$, then
$\lambda_{w,\sigma}$ permutes the $|w|+1$st nonzero letter and the
following letter by $\sigma$. Otherwise $\lambda_{w,\sigma}$ does nothing. ( e.g. $\lambda_{21,(01)\wr (12)}(\cdots 001020010) = \cdots 012020010$ and $\lambda_{21,(01)\wr (12)}(\cdots 002010010) = \cdots 002010010$)

For a word $w$ of length $k$, define the group $\L^w_{m,k}$
generated by $\lambda_{w,\sigma}$ as $\sigma$ ranges over
$\Sym(m)\wr\Sym(m-1)$. Define the group $\L_{k,m}$ to be the group
generated by the $\L^w_{k,m}$ for all words $w$ of length
$k$. Define further $\L_{-1,m}=\Sym(m)$.

Later, we will consider random walks on the mother group whose step
distribution is a convex combination of uniform measures on the subgroups
$\L_{k,m}$ for various $k$'s.

\begin{lemma}[\cite{AAV}]
  For each $w$, $\L^w_{k,m} \approx \Sym(m)\wr\Sym(m-1)$. The
  group $\L_{k,m}$ is a subgroup of
  $\M_{k,m}$ and is the direct product of $\L^w_{k,m}$ for
  $w\in\{1,\dots,m-1\}^k$. Moreover, the mother group $\M_{k,m}$ is generated by
the subgroups $\{\L_{m,\ell}\}_{\ell\leq k}.$
\end{lemma}

\begin{definition}
The {\bf fractal mother graphs} $\G(d,m,n)$ are the
Schreier graphs of the mother group $\M_{d,m}$ acting on
level $n$ of $\tree _m$ with the generating set
$\bigcup_{k\leq d} \L_{k,m}$ (counting multiplicity).
These are regular graphs with vertex degree
$$
deg(v)= \sum_{k=-1}^d | \L_{k,m}| = m! + \sum_{k=0}^d (m-1)^k
(m-1)!(m!)^{m-1}
$$

$\G(d,m,\infty)$ is the connected component of the zero ray
of the Schreier graph of the action of $\M_{d,m}$ on the
boundary. The figures in this paper depict various
instances of $\G(d,m,n)$.
\end{definition}

\section{Basic resistance properties of fractal mother graphs}

The vertices of $\G(d,m,n)$ are the set of words in
$\{0,\ldots,m-1\}^n$. There is a natural embedding $\G(d,m,1)\subset \G(d,m,2)\cdots\subset
\G(d,m,\infty)$ of the vertex sets of these graphs, given by
adding zeros to the beginning of words. This embedding
preserves the graph structure in the sense that the graph
induced by the vertices of $\G(d,m,i)$ in $\G(d,m,{i+1})$ is isomorphic
to $\G(d,m,i)$ with some loops erased.


We introduce some notation: The $n$-digit vertices $o=00\ldots 00$ and $\oo_n(x)=x0\ldots 00$ ($x\in \{1,2,\cdot,m-1\}$) are  called the \textbf{root} and $x$-\textbf{antiroot} of $\G(d,m,n)$. $\oo_n(1)$ is called simply the \textbf{antiroot} of $\G(d,m,n)$ and denoted by $\oo_n$. For
$w\in \G(d,m,n)$ let $\#w$ denote the number of nonzero letters in $w$.
Finally, let $\res_{d,m,n}(a,b)$ denote the resistance in $\G(d,m,n)
$ between vertices or vertex sets $a,b$.

\begin{proposition}\label{p.res}

\ \\  \vspace{-2em}

\begin{enumerate}[(a)]
      \item \label{p.res.c} For $x\in \G(d,m,n)$ we have
          $\res_{d,m,n}(o,x) \le \sum_{s:x_s\not=0}
          \res_{d,m,s}(o,\oo_s)$.
      \item \label{p.res.d} For $x,y\in \G(d,m,n)$ we have
          $\res_{d,m,n}(x,y) \le \sum_{s=1}^n
          (\one(x_s\not=0)+\one(y_s\not=0))
          \res_{d,m,s}(o,\oo_s)$.
\end{enumerate}
\end{proposition}

\begin{proof} \ \\
\noindent \eqref{p.res.c}   Let $x=x_n\ldots x_1$, and
$z_s=x_n\ldots x_s$ followed by $s-1$ zeros. The triangle
inequality for resistances gives
$$
\res_{d,m,n}(o,x)\le \sum_{s=1}^n \res_{d,m,n}(z_{s-1},z_s).
$$
The terms where $x_s=0$ are 0.  To complete the proof, we
need to show that
$$
\res_{d,m,n}(z_{s-1},z_s) \le \res_{d,m,s} (o,\oo_s).
$$
This is because the map (mapping words of length $s$ to
words of length $n$) that sends a word $v$ to $x_n\ldots
x_{k-1} v$ is an injection of vertices and extends to an
injection of edges (apart from loops). Moreover,
$$o\mapsto z_s \quad \mbox{
 and }\quad \oo_\ell(x_s) \mapsto z_{s-1}.$$
 Thus by Rayleigh's monotonicity we get
$$
\res_{d,m,n}(z_{s-1},z_s) \le \res_{d,m,s} (o,\oo_s(x_s))
$$
and the right-hand side equals $\res_{d,m,s} (o,\oo_s)$ by
symmetry.

\noindent \eqref{p.res.d} This follows from part
\eqref{p.res.c} and the triangle inequality for
resistances.
\end{proof}
Let
$$\rr(d,m,n)=\max_{x,y\in \G(d,m,n)}\res_{d,m,n}(x,y)
$$
\begin{proposition}\label{p.res2}
\ \\  \vspace{-2em}
\begin{enumerate}[(a)]
      \item\label{p.res2.a} For $n'\le n$ and $x,y<m^{n'}$ we have
      $\res_{d,m,n}(x,y)\le
      \res_{d,m,n'}(x,y)$;
      \item \label{p.res2.b}
      There exists a $c$ depending on $m,d$ only so that
      for $k\ge 1$ we have
      $\rr(d,m,n+k)\le
      c^k \rr(d,m,n)$.
\end{enumerate}
\end{proposition}

\begin{proof} \ \\
\noindent \eqref{p.res2.a} This follows from Rayley's
monotonicity principle and the fact that, apart from loops,
$\G(d,m,n')$ is a subgraph of $\G(d,m,n)$.

%
\smallskip
\noindent \eqref{p.res2.b}  It suffices to prove this for
$k=1$, the rest follow by induction. Consider the vertices
$x <m^n$. Fix a letter $b$. Note that $x\mapsto bx$ is a
graph homomorphism (except for loops). In particular, by
Rayleigh's monotonicity principle for $x,y<m^n$ we have
\begin{equation}\label{e:resi}
\res_{d,m,n+1}(bx,by)\le \res_{d,m,n}(x,y).
\end{equation}
Recall that $\oo_n$ is the antiroot $100\ldots 000$ in
$\G(d,m,n)$. For letters $a,b$ the triangle inequality for
resistances implies
$$
\res_{d,m,n+1}(ax,by)\le
\res_{d,m,n+1}(ax,a\oo_n)+\res_{d,m,n+1}(a\oo_n,b\oo_n)+\res_{d,m,n+1}(b\oo_n,by)
$$
There is an edge between $a\oo_n$ and $b\oo_n$. So by
\eqref{e:resi} we have the upper bound
\[
\res_{d,m,n+1}(ax,by)\le
\res_{d,m,n}(x,\oo_n)+1+\res_{d,m,n}(y,\oo_n)\le (2+c)
\rr(d,m,n). \qedhere
\]
\end{proof}

\section{Flow construction}

The goal of this section is to give an upper bound on the
maximal resistance in $\G(d,m,n)$. We include the proof of
the following simple fact for completeness.

\begin{lemma}\label{l:zero_res}
We have  $\rr(0,m,n)\le 2^n-1$.
\end{lemma}

\begin{proof}
It suffice to prove that there is a path of length $2^n-1$
to any $n$-digit number $\ell$. By the symmetry of the
mother group, it suffices to show this for numbers
containing the digits $0,1$ only; in particular, it
suffices to show this for the $m=2$ mother group. For this,
for the binary string $x=x_n\ldots x_1$ consider the
inverse Gray code representation $y=y_n\ldots y_1$ with
bits $y_k=x_n+\ldots +x_k$ . By checking how the two
generators act on the $y$, one sees that graph indexed by
the  $y$'s is the simple path $0,1,2,\ldots 2^n-1$ (with loops attached at the endpoints).
\end{proof}

\begin{proposition}\label{p:degree0bound}
We have  $c\left(\frac{m}{m-1}\right)^n \le \rr(0,m,n)\le
c' \left(\frac{m}{m-1}\right)^n$, where $c,c'$ depend on
$m$ only.
\end{proposition}

\begin{figure}[t]
\begin{center}
\includegraphics*[width=300pt]{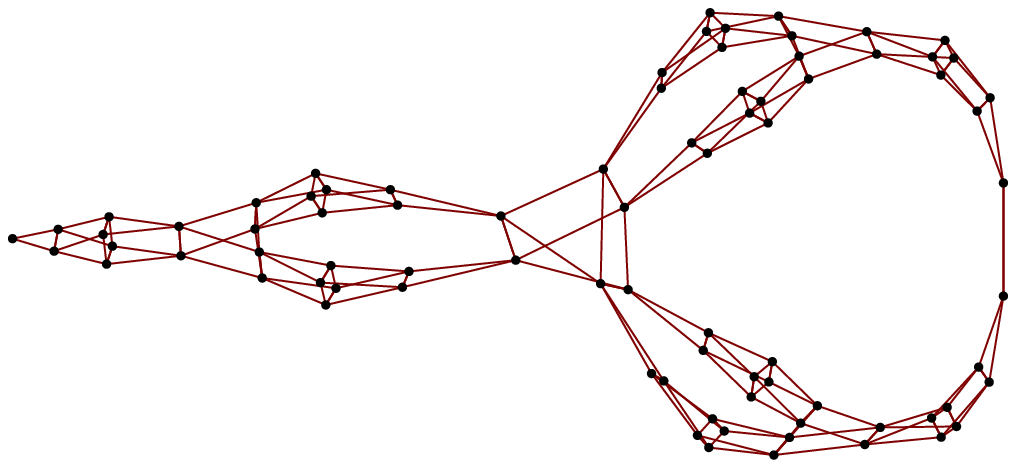}\\
The fractal mother graph $\G(0,3,4)$
\end{center}
\end{figure}

\begin{proof}
First let $d$ be general, and the current flow from $o$ to
the set of antiroots $\bar{O}_n=\{\bar{o}_n(1),\ldots
\bar{o}_n(m-1)\}$ in $\G(d,m,n)$. Note that the graph, $o$
and $\bar O_n$ are symmetric under the action of
$\Sym(m-1)^n$ which changes the nonzero letters of the
strings of each vertex. Since the current flow is unique,
it must be invariant under this action.  In particular, if
we identify vertices of the graph that are in the same
orbit, the resistances will not change (normally, they
could only decrease by Raleigh's monotonicity) After this
identification, we get a weighted version of $\G(d,2,n)$
with vertices given by binary strings. Each vertex $y$ is
the image of $(m-1)^{\#y}$ vertices in $\G(d,2,n)$. For
$d=0$ this is graph is a simple path (Lemma
\ref{l:zero_res}), and the resistance between the endpoint
can be computed explicitly. Here we estimate it by a quick
argument: each vertex has one or two incident edges. The
weight of these edges is within constants of the weight of
the vertex. By the series formula, the resistance is given
by the sum over the edges of $w(e)^{-1}$, which is
therefore within constants of the sum over the vertices of
$w(v)^{-1}$, that is
$$
\sum_{v=0}^{2^m-1}(m-1)^{-\#v}=\left(\frac{m}{m-1}\right)^n,
$$
computed via the binomial formula. Finally, note that symmetry implies
$$\res_{d,m,n}(o,\bar O_n) \le \res_{d,m,n}(o,\bar o_n)\le (m-1)^2 \res_{d,m,n}(o,\bar O_n).
$$
for the last inequality, note that the current flow to
$O_n$ has strength $1/(m-1)$ going into $\bar o_n$, so the
energy of $(m-1)$ times this flow is an upper bound for the
resistance between $o$ and $\bar o_n$. Finally, Proposition
\ref{p.res} \eqref{p.res.d} implies the bound on $\bar r$.
\end{proof}

One way to get a transient graph is to construct a product
structure; for example $\mathbb Z \times \mathbb Z \times
\mathbb Z$ is transient. We would like to construct
something that resembles the product of the Schreier graphs
of a degree $d$ group and a degree $d'$ group. It will turn
out to be a subgraph of the zero-ray Schreier graph of the
mother group of degree $d$ for greater $n$.

Recall that  $\rr(d,m,n)$ denotes the maximal resistance
between vertices in $\G(d,m,n)$. The core of the argument
is contained in the following proposition.

\begin{proposition}\label{p:flow_bound} For any nondecreasing sequence
$\gamma_s\ge 1$ and $n\ge 3$ we have
\begin{equation}\label{rmdn} \rr(d,m,n)\le  c
\sum_{s=1}^{n-2}\left[\frac{\rr(d',m,\lfloor
\gamma_{s}\rfloor)}{s^{d-d'-1}}+ \frac{\rr(d,m, s)
}{m^{\gamma_s}}\right],
\end{equation}
where $0\le d' \le d-1$ is an integer and $c$ is a constant
depending on $m$ and $d$ only.
\end{proposition}

\begin{proof}
We first assume that the $\gamma_s$ are integers,
$\gamma_1=\ldots =\gamma_d=1$, and $n\ge d+2$.
Let $a,a'$ be two vertices in $\G(d,m,n)$. Thompson's principle says that
resistance between two vertices equals the energy of the
minimal energy unit flow between them.
We will construct a unit flow from $a$ to $a'$  which has
the right bounds. We will construct this flow as a sum of
many parts.

Let $a_k$ denote $a$ with its lowest $k$ digits erased. Let
$a^*_k$ denote $a_k$ with its lowest digit permuted
cyclically by $(0\ldots m-1)$. Let $\sigma$ be the smallest
number so that $\gamma_\sigma+ 1+\sigma \ge n$; we note
that $\sigma \ge d$.

Recall that  $\#y$ denotes the number of nonzero digits of
$y$. For $d-1\leq s<\sigma$ let
$$X_s=\{a^*_{\gamma_s+s+1}x\,:\,0\le x<m^{\gamma_s}\}$$
where $a^*_{\gamma_s+s+1}x$ denotes the concatenation of
$a^*_{\gamma_s+s+1}$ and the $\gamma_s$-digit version of
$x$.
Let
$$Y_s=\{0\le y< m^s\,:\,\#y= d-d'-1\}.$$ We will think of numbers in $Y_s$ as written in their
base-$m$, $s$-digit form.
Finally, define
$$X_\sigma=\{x :\,0\le x<m^{n-\sigma}\}$$
For $x\in X_s$ and $y\in Y_s$ we denote $x1y$ the obvious
concatenation, and let $X_t1Y_s$ denote the set of such
concatenations; let $X_t01Y_s$ be defined similarly.

We now consider flows $$a\to X_{d-1}1Y_{d-1} \to
X_{d}01Y_{d-1}\to X_{d}1Y_d\to X_{d+1}01Y_d\to
X_{d+1}1Y_{d+1}\to\cdots\to X_{\sigma}01Y_{\sigma-1}.
$$ each transporting uniform measure from one set to
uniform measure on the next. These vertex sets and flows
will not be disjoint, but we will be able to control
overlaps.

Let $\phi$ be the sum of these flows. Let $\phi'$ be the
sum of the same flows constructed for $a'$, and let $\bar
\phi$ be its reversal. Then $\phi+\bar \phi$ is a unit flow
from $a$ to $a'$.

Let $\energy(\phi)$ denote the energy of the flow $\phi$,
i.e. $\sum_e \phi(e)^2$. To bound the energy
of $\phi$, we first bound the energy of its components;
these alternate between two types.

\begin{description}
\item[Horizontal spread, $X_{s-1}1Y_{s-1}$ to
$X_{s}01Y_{s-1}$.] Fix $y\in Y_{s-1}$. Consider edges of
the form $(a^*_{s+\gamma_s+2}x1y,a^*_{s+\gamma_s+2}(g.x)1y)$
where $g$ is a generator of degree at most $d'$ and
$x,g.x<m^{\gamma_{s}+2}$. The subgraph of $\G(d,m,n)$
spanned by these edges is isomorphic (up to loops) to
$\G(d',m,\gamma_{s}+2)$ through the map $x\mapsto
a^*_{s+\gamma_s+2}x1y$.

We call such edges {\bf horizontal}. By Thompson's
principle for any $x_1\in X_{s-1}$ and $x_2\in X_{s}$ there
exists a unit flow $\phi^h_{s,x_1,x_2,y}$ from $x_11y$ to
$x_201y$ to along these edges satisfying
$$
\energy(\phi^h_{s,x_1,x_2,y})\le \rr(d',m,\gamma_s+2).
$$
Let $\phi^h_{s,y}$ denote the average of these flows over
all $x_1\in X_{s-1},\,x_2\in X_{s}$. By convexity of
energy, we have
$$
\energy(\phi^h_{s,y})\le \rr(d',m,\gamma_s+2).
$$
We now average these flows over all $y\in Y_{s-1}$ to get
the flow $\phi^h_s$. Since these flows have disjoint
support, and since $|Y_{s-1}|\ge c s^{d-d'-1}$ we get that
\begin{eqnarray*}
\energy(\phi^h_{s})&\le& c s^{1+d'-d} \rr(d',m,\gamma_s+2)
\le c s^{1+d'-d} \rr(d',m,\gamma_{s}).
\end{eqnarray*} The last inequality follows from
Proposition \ref{p.res2.a}. The flow $\phi^h_s$
transports uniform measure form $X_{s-1}1Y_{s-1}$ to
uniform measure on $X_{s}01Y_{s-1}$.

A trivial modification of the argument gives a flow
$\phi_\sigma$ from  $X_{\sigma-1}1Y_{\sigma-1}$ to
$X_{\sigma}01Y_{\sigma-1}$ with energy bounded by
\begin{eqnarray*} \energy(\phi^h_{s}) \le c s^{1+d'-d}
\rr(d',m,\gamma_{\sigma}).
\end{eqnarray*}

\item[Vertical spread, $X_{s}01Y_{s-1}$ to $X_{s}1Y_{s}$:]
Fix  $x\in X_{s}$. Consider edges of the form $(xy_s,xg.y)$
where $g$ is a generator of degree at most $d$, and
$y,g.y<m^{s+1}$ and $y$,$g.y$ above are in their
$s+1$-digit form. The subgraph of $\G(d,m,n)$ spanned by
these edges is isomorphic (up to loops) to $\G(d,m,s+1)$.
We call such edges {\bf vertical}. By Thompson's principle
for any $y_1\in Y_{s-1}$ and $y_2\in Y_{s}$ there exists a
unit flow $\phi^v_{s,x,y_1,y_2}$ from $x01y_1$ to $x1y_2$
along these edges satisfying
$$
\energy(\phi^v_{s,x,y_1,y_2})=\res_{d,m,s+1}(01y_1,1y_2).
$$
Let $\phi^v_{s,x}$ denote the average of these flows over
all $y_1\in Y_{s-1}$ and $y_2\in Y_{s}$. By convexity of
energy, we have
$$
\energy(\phi^v_{s,x})\le  \rr(d,m,s+1).
$$
We now average these flows over all $x\in X_{s}$ to get the
flow $\phi^v_s$. Since these flows have disjoint support,
we get that
$$
\energy(\phi^h_{s})\le \frac{1}{m^{\gamma_{s}}}
\rr(d,m,s+1)\le \frac{c}{m^{\gamma_{s}}} \rr(d,m,s).
$$
The last inequality follows from Proposition \ref{p.res2.a}. The flow $\phi^v_s$ transports uniform measure on
$X_{s}01Y_{s-1}$ to uniform measure on $X_{s}1Y_{s}$.
\end{description}

We consider the flow from $o$ to $V_{4,s}$ constructed piecewise
as above, namely
\begin{equation} \label{phi}
  \phi = \phi_a + \phi_{d}^h+\phi_{d}^v +\ldots +
  \phi_{\sigma-1}^v
+\phi_\sigma^h.
\end{equation}
Here $\phi_a$ a unit current flow from $a$ to uniform
measure on $X_{d-1}1Y_{d-1}$, and its energy is bounded by
some constant.

We first note that the flows $\varphi^h_s$ for different
$s<\sigma$ are vertex-disjoint. This is because these flows
move between vertices of the form $x1y$, with $y\in
Y_{s-1}$ fixed. As long as one can determine $s$ from
$x1y$, the vertex-disjointness follows. But since
$\#y=d-d'-1$, the $d'+2$nd nonzero digit  from the right in
$x1y$ is in position $s$.

Similarly, we show that the flows $\varphi^v_s$ are vertex
disjoint for different $s<\sigma$. Again, these flows are
vertex-disjoint. They move between vertices of the form
$xy$ with $y<m^{s+1}$ is in $s+1$-digit form, and $x\in
X_s$ is fixed. This is why we permuted a bit in $a^*$: the
highest digit in $xy$ that differs from the digit at the same
position in $a$ is exactly at position $s+\gamma_s+2$.
Since the map $s\mapsto s+\gamma_s$ is strictly increasing,
the value of $s+\gamma_s$ determines the value of $s$.

Note that the energy of the sum of flows $\sum \phi_i$ is
bounded above by $\sum m_i\energy(\phi_i)$ where $m_i$ is
the number of flows (including itself) that the flow
$\phi_i$ shares an edge with. Breaking up the flow
$\varphi$ into four parts, namely the sum of the horizontal
flow terms, the vertical terms, and $\phi_*$ as well as
$\phi_\sigma$
$$
\energy(\varphi)\le 4\energy (\varphi_0) +4 c
\sum_{s=d}^{\sigma}\frac{\rr(d',m,
\gamma_{s})}{s^{d-d'-1}}+ 4c\sum_{s=d}^{\sigma-1}
\frac{\rr(d,m, s) }{m^{\gamma_s}}
$$
Note that the flow $\varphi_0$ can be chosen so that its
energy is bounded by a constant depending on $d,m$ only. We
can remove the additive constants by increasing the
multiplicative one and using that $r(d',m,\cdot)$ is
bounded below. Then Thompson's principle gives
$$\res_{d,m,n}(a,a')\le \energy(\phi+\overline\phi) \le
2\energy(\phi)+2\energy(\overline\phi)\le c
\sum_{s=d}^{\sigma}\frac{\rr(d',m,
\gamma_{s})}{s^{d-d'-1}}+ c\sum_{s=d}^{\sigma-1}
\frac{\rr(d,m, s) }{m^{\gamma_s}}.
$$
Since $\sigma\le n-2$, the inequality \eqref{rmdn} follows
for $\gamma_s$, $n$ satisfying the assumptions at the
beginning of the proof. Using the fact that
$\rr(d',m,\cdot)$ is bounded below, and Proposition
\ref{p.res2} \eqref{p.res2.b} it is straightforward to see
that \eqref{rmdn} holds with a modified constant $c$ for
general $\gamma_s$, $n$.
\end{proof}

The proof of the following standard fact is a simple
exercise.

\begin{fact}\label{f:transience}
Let $G$ be a bounded degree graph with a vertex $o$ with an
infinite sequence of vertices whose resistance to $o$ is
bounded. Then $G$ is transient.
\end{fact}

\begin{proposition}\label{p:mainresbound}For $d\ge 1$, $m\ge 2$ we have
$$\rr(d,m,n)\le \begin{cases}
 c &\text{for } \{d\ge 4,m\ge 2\}  \text{ or }\{d\ge 3, m \ge 3\},\\
     c\log^2 n &\text{for } \{d=3,m=2\},\\
cs^{3-d-\log_m (m-1)} &\text{for } \{1\le d\le 2\}.\\
     \end{cases}
$$
Where $c$ depends on $m,d$ only. In particular, in the
first case, $\G(d,m,\infty)$ is transient.
\end{proposition}
\begin{proof}
Set $d'=0$, and set $$r(s)=r(d,m,s) = \max_{1\le \ell \le s}\rr(d,m,\ell)$$
to guarantee that $r$ is nondecreasing. Proposition
\ref{p:flow_bound} together with Lemma \ref{l:zero_res}
gives
$$
r(n+2) \le c\sum_{s=1}^n
\frac{\left(\frac{m}{m-1}\right)^{\gamma_s}}{s^{d-1}} +
\frac{r(s)}{m^{\gamma_s}}.
$$
Since we are free to choose $\gamma_s$, we will find the
approximate minimizer of the right hand side. When the two
terms are equal, we are off from the minimum by at most a
factor of two. The calculation gives the choice
$$
\gamma_s=\log_{m^2/(m-1)}(r(s)s^{d-1})
$$
which is indeed nondecreasing and gives
$$
r(n+2) \le c\sum_{s=1}^n \frac{r(s)^{1-\eta}}
{s^{(d-1)\eta}},\qquad \eta = (2-\log_{m}(m-1))^{-1}.
$$
let $f(s)$ denote the summand on the right hand side. We
extend $r$ to a piecewise linear continuous function on
$[1,\infty)$. For this function we have
$$
r(t)\le r(\lceil t \rceil) \le  \sum_{s=1}^ {\lceil t
\rceil-2} f(s)\,ds  \le \int_1^{\lceil t \rceil-1} f(s)\,ds
\le \int_{1}^{t} f(s)\,ds
$$
for $t\ge 3$. In particular, adding a constant, we can make
the inequality valid for all $t\ge 1$:
\begin{equation}\label{e:integralr}
r(t)\le c+c\int_{1}^{t}\frac{r(s)^{1-\eta}}
{s^{(d-1)\eta}}\,ds.
\end{equation}
consider the solution of the differential equation
$$
u'(s) =\kappa \frac{u(s)^{1-\eta}} {s^{(d-1)\eta}}
$$
given by
$$
u(s)=\begin{cases}
 {\left(\frac{\kappa \eta(1- s^
        {1 - (d-1)\eta })}{(d-1)\eta-1 }
     \right) }^{\frac{1}{\eta }} &\text{for } (d-1)\eta>1,\\
     {\left( \eta \,\kappa \,\log (s) \right) }^
  {\frac{1}{\eta }} &\text{for } (d-1)\eta=1,\\
{\left(\frac{{{\kappa \eta s}}^
        {1 - (d-1)\eta }}{1 - (d-1)\eta }
     \right) }^{\frac{1}{\eta }} &\text{for } (d-1)\eta<1.\\
     \end{cases}
     $$
for $\kappa$ large enough (and greater than the constant in
front of the integrand in \eqref{e:integralr}) we have
$u(2)$ will be greater than the right hand side of
\eqref{e:integralr} evaluated at $s=2$. For such $\kappa$
we claim that $u(t)>r(t)$ for all $t\ge 2$. Assume the
contrary;  since $r,u$ are continuous, then there is a
smallest $t\ge 2$ so that $u(t)=r(t)$. Therefore
$$
r(t) \le c + c\int_1^t \frac{r(s)^{1-\eta}}
{s^{(d-1)\eta}}\,ds < u(2) + \kappa \int_2^t
\frac{u(s)^{1-\eta}} {s^{(d-1)\eta}}\,ds=u(t),
$$
a contradiction. In terms of $d$ and $m$, the inequality
deciding the cases above is $d\ge 3-\log_m (m-1)$. The
claim of the corollary follows.
\end{proof}

\section{From transience to nontrivial Poisson boundary}

Let $K$ denote the subgroup of $\M_{d,m}$ generated by all
generators except those of degree $-1$. Note that all
elements of $K$ fix the zero ray. Let $H$ denote the
subgroup of $K$ generated by all generators excepts those
of degree $-1$ or degree $d$.

\begin{proposition}
For every $g$ in the mother group, as $n\to \infty$ the
$n$th level section at the zero ray is eventually constant
and is in $K$. We call this germ$(g)$.
\end{proposition}

Note that $germ(k)=k$ for $k\in K$.

\begin{proof}
We prove this by induction on $|g|$. Write
$$
g=\ang{g_0=x,\ldots ,g_{m^n-1}} \tau
$$
the decomposition of $g$ at level $n$. Assume that the
claim holds for $g$, i.e. $x\in K$; we will show it for
$gs$ where $s$ is a generator. Since $x\in K$, we have
$0.x=0$ and the first-level section satisfies $x(0)=x$.
Then for a generator $s=\ang{s_0,\ldots s_{m^n-1}}\sigma$
we have
$$
gs=\ang{x y,\ldots }\tau\sigma
$$
where $y=s_{0.\tau}$ is itself a generator or the identity.
At level $n+1$ we have the expansion
$$gs=\ang{xy(0),\ldots}\tau'\sigma'.$$
 Now for every generator $y$ if $y\notin K$, then
$y(0)=id$. In this case we have $xy(0)=x(0)=x$, and the
level $n+1$ expansion is of the form
$gs=\ang{x,\ldots}\tau'\sigma'$. Note that $x(0)=x$ implies
that the level sections of $g$ at the zero ray remain $x$
from level $n+1$ on.

Otherwise, if $y\in K$,  note that $0.x=0$, and $y(0)=y$,
and therefore $xy(0)=x(0)y(0.x)=xy$, so the level sections
of $g$ remain $xy$ from level $n$ on, the completing the
proof.
\end{proof}

We define lamp$(g)$ as the quotient of germ$(g)$ by the
subgroup group generated by all generators of degree at
most $d-1$. The function lamp$()$ takes values in the coset
space $K/H$.

\begin{proposition}[Lighting a single lamp]
Let $g\in M$ and let $s$ be a generator. Then we have
$$
\lamp(gs)=
\begin{cases}
\lamp(g)\lamp(s) & \mbox{ if  $g$ fixes the zero ray}\\
\lamp(g) & \mbox{ otherwise.}
\end{cases}
$$
\end{proposition}

\begin{proof}
We use the notation of the previous proof. The proof implies
\begin{equation}\label{gy}
\germ(gs)=\germ(g)\germ(y)
\end{equation}
(for the generator $y$, if  $y\notin K$ then
$\germ(y)=id$).

If  $g$ fixes the zero ray, then $0.\tau=0$ and $y=s$, so
the claim follows from \eqref{gy}.

 If $g$ does not fix
the zero ray, then at a high enough level
$0.\tau\not=\tau$. Working on a level where this happens
and also the zero-ray section of $g$ is stable, we get
$y=s_{0.\tau}$, and therefore $y$ is of lower degree than
$s$. Thus $\lamp(y)=idH$ (trivial), and the claim follows
from \eqref{gy}.
\end{proof}

We are now ready to prove theorem \ref{t:liouville-tran},
which states that if the Schreier graph of the mother group
the zero ray is transient, then no non-degenerate random
walk with finitely supported step distribution has the
Liouville property.
\begin{proof}[Proof of Theorem \ref{t:liouville-tran}]
Note that transience is quasi-isometry invariant, and
changing generators is a quasi-isometry of the Schreier
graph. Thus the Schreier graph for the random walk
$Y_1\cdots Y_n$ in question is also transient. Thus
$Y_1\cdots Y_n$ in question fixes the zero ray only
finitely many times. In particular, lamp$(Y_1\cdots Y_n)$
stabilizes at some random value $L$ with probability 1.

We now claim that $L$ takes all possible value in $K/H$
with positive probability, and hence the Poisson boundary
is nontrivial. Indeed, let $l\in K/H$ with $|l|=n$. Then
with positive probability $Y_1\cdots Y_n=l$. Moreover, by
transience of the Schreier graph, the independent walk
$Y_{n+1}\cdots Y_{n+k}$, $k\ge 0$ never fixes the zero ray
for $k \ge 1$ with positive probability. In this case
$L=\ell$.

We have shown that the Poisson boundary is nontrivial. Now,
by the standard construction, function $x\mapsto
P_x(L=id)$, where $x$ is the starting point of the random
walk, is bounded, harmonic and nonconstant.
\end{proof}
%
%
%
%
\paragraph{Acknowledgements.} Part of this work was done while the
first author was at the University of Toronto. The work of
the second author was supported by the Canada Research Chair program and
the NSERC Discovery Accelerator program. In addition, he
thanks the Department of Stochastics, Technical University
of Budapest where this work was completed.
\bibliographystyle{dcu}
\bibliography{amenability}

\bigskip\bigskip\bigskip\noindent
\begin{minipage}{0.49\linewidth}
Gideon Amir
\\Department of Mathematics
\\Bar Ilan University
\\Ramat Gan, 52900 Israel
\\{\tt amirgi@macs.biu.ac.il }
\end{minipage}
\begin{minipage}{0.49\linewidth}
B\'alint Vir\'ag
\\Departments of Mathematics and Statistics
\\University of Toronto
\\Toronto ON~~M5S 2E4, Canada
\\{\tt balint@math.toronto.edu}
\end{minipage}

\end{document}